\documentclass[a4paper,11pt]{amsart}
\usepackage{graphicx}
\usepackage{mathptmx}
\usepackage{amsmath}
\usepackage{amssymb}
\usepackage{enumitem}
\usepackage{xcolor}

\newmuskip\pFqmuskip

\newcommand*\pFq[6][8]{%
  \begingroup % only local assignments
  \pFqmuskip=#1mu\relax
  \mathcode`=\string"8000
  \begingroup\lccode`\~=`\,
  \lowercase{\endgroup\let~}\pFqcomma
  F^{#2}_{#3}{\left(\genfrac..{0pt}{}{#4}{#5}\bigg|#6\right)}%
  \endgroup
}
\newcommand{\pFqcomma}{\mskip\pFqmuskip}

\newtheorem{theorem}{Theorem}

\begin{document}

\title[]{Degenerate $r$-Bell polynomials \\ arising from degenerate normal ordering }
\author{Taekyun  Kim$^{1,*}$}
\address{Department of Mathematics, Kwangwoon University, Seoul 139-701, Republic of Korea}
\email{kwangwoonmath@kw.ac.kr}

\author{Dae San  Kim$^{2,*}$}
\address{Department of Mathematics, Sogang University, Seoul 121-742, Republic of Korea}
\email{dskim@sogang.ac.kr}

\author{Hye Kyung  Kim$^{3,*}$}
\address{Department of Mathematics Education, Daegu Catholic University, Gyeongsan 38430, Republic of Korea}
\email{hkkim@cu.ac.kr}

\thanks{* are corresponding authors}

\begin{abstract}
Recently, Kim-Kim introduced the degenerate $r$-Bell polynomials and investigated some results which are derived from umbral calculus. The aim of this paper is to study some properties of the degenerate $r$-Bell polynomials and numbers via boson operators. In particular, we obtain two expressions for the generating function of the degenerate $r$-Bell polynomials in $|z|^{2}$, and a recurrence relation and Dobinski-like formula for the degenerate $r$-Bell numbers. These are derived from the degenerate normal ordering of a degenerate integral power of the number operator in terms of boson operators where the degenerate $r$-Stirling numbers of the second kind appear as the coefficients.
\end{abstract}

\maketitle

\markboth{\centerline{\scriptsize Degenerate $r$-Bell polynomials arising from degenerate normal ordering}}
{\centerline{\scriptsize  T. Kim, D. S. Kim and  H. K. Kim}}

\section{Introduction and preliminaries}
It turns out that it is fascinating and fruitful to study varioius degenerate versions of some special polynomials and numbers (see [6-12] and the references therein), which has its origin in the work of Carlitz in [2]. They have been explored by using various different methods and led to the introduction of degenerate gamma functions and degenerate umbral calculus. \par
The aim of this paper is to study some properties of the degenerate $r$-Bell polynomials and numbers via boson operators. In particular, we obtain two expressions for the generating function of the degenerate $r$-Bell polynomials in $|z|^{2}$, and a recurrence relation and Dobinski-like formula for the degenerate $r$-Bell numbers. These are derived from the degenerate normal orderings of a degenerate integral power of the number operator in terms of boson operators where the degenerate $r$-Stirling numbers of the second kind appear as the coefficients.\par
In more detail, the outline of this paper is as follows. In Section 1, we recall the degenerate exponentials, the degenerate $r$-Stirling numbers of the second kind and the degenerate $r$-Bell polynomials. We remind the reader of the boson operators, the number operators, the normal ordering of an integral power of the number operator in terms of boson operators and the degenerate normal ordering of a degenerate integral power of the number operator in terms of boson operators. Section 2 is the main result of this paper. In Theorem 1, we state that the degenerate normal orderings of the `degenerate integral powers of the number operator' $(\hat{n}+r)_{n,\lambda}$ and $(\hat{n})_{n-r,\lambda}(a^\dag)^{r}a^{r}$ in terms of the boson operators with the coefficients given by the degenerate $r$-Stirling numbers of the second kind. Let $f(t)=\langle{z|e_{\lambda}^{r+\hat{n}}(t)|z \rangle}$. We show that $\langle{z|(\hat{n}+r)_{n,\lambda}|z \rangle}$ is equal to $\phi_{k,\lambda}^{(r)}(|z|^2)$ in Theorem 2 and hence that $f(t)$ is the generating function of that in Theorem 3, where $\phi_{k,\lambda}^{(r)}(x)$ is the degenerate $r$-Bell polynomial. We derive a differential equation for $f(t)$ in Theorem 4 and thereby get another expression for $f(t)$ in Theorem 5, which in turn yields an explicit generating function of $\phi_{k,\lambda}^{(r)}(|z|^2)$. In Theorem 6, we obtain a recurrence relation for the degenerate $r$-Bell numbers $\phi_{n,\lambda}^{(r)}=\phi_{n,\lambda}^{(r)}(1)$. Finally, we get a Dobinski-like formula for the degenerate $r$-Bell numbers from the representation of the coherent state in terms of the number operator in Theorem 7. For the rest of this section, we recall the facts that are needed throughout this paper. \par
For any nonzero $\lambda\in\mathbb{R}$, the degenerate exponentials are defined by
\begin{equation}\label{eq01}
\begin{split}
e_\lambda^x(t)=(1+\lambda t)^{\frac{x}{\lambda}}=\sum_{k=0}^\infty (x)_{k,\lambda}\frac{t^k}{k!},\quad e_\lambda(t)=e_\lambda^1(t),
\end{split}
\end{equation}
where
\begin{equation}\label{eq02}
\begin{split}
(x)_{0,\lambda}=1, \ \ (x)_{k,\lambda}=x(x-\lambda)\cdots(x-(k-1)\lambda),\quad (k\geq1),\quad \ \ (\text{see \ [2, 6-9,12]}).
\end{split}
\end{equation}

For $r\in\mathbb{Z}$ with $r\geq0$, the degenerate $r$-Stirling numbers of the second kind are defined by Kim-Kim as

\begin{equation}\label{eq03}
\begin{split}
(x+r)_{n,\lambda}=\sum_{k=0}^n \begin{Bmatrix}n+r \\ k+r \end{Bmatrix}_{r,\lambda} (x)_k, \quad (n\geq0), \quad \ \ (\text{see  \ [10-12]}),
\end{split}
\end{equation}
where $(x)_0=1, (x)_n=x(x-1) \cdots (x-n+1),\quad (n \ge 1)$.

From \eqref{eq03}, we note that

\begin{equation}\label{eq04}
\begin{split}
\frac{1}{k!}(e_\lambda (t)-1)^k e_{\lambda}^r (t)=\sum_{n=k}^\infty \begin{Bmatrix}n+r \\ k+r \end{Bmatrix}_{r,\lambda} \frac{t^n}{n!},\quad (k\geq0), \quad \ \ (\text{see \ \cite{11,12}}).
\end{split}
\end{equation}

Note that $\lim_{\lambda\rightarrow 0} \begin{Bmatrix}n+r \\ k+r \end{Bmatrix}_{r,\lambda}=\begin{Bmatrix}n+r \\ k+r \end{Bmatrix}_{r}$ are the $r$-Stirling numbers of the second kind.
It is well known that the $r$-Stirling number of the second kind $\begin{Bmatrix}n+r \\ k+r \end{Bmatrix}_{r}$ counts the number of partitions of the set $[n]=\{1,2,\cdots,n\}$ into $k$ nonempty disjoint subsets in such a way that the numbers $1,2,3,\cdots,r $ are in distinct subsets.

The degenerate $r$-Bell polynomials $\phi_{n,\lambda}^{(r)}(x)$ are defined by
\begin{equation}\label{eq05}
\begin{split}
e_\lambda^r (t)e^{x(e_\lambda(t)-1)}=\sum_{n=0}^\infty \phi_{n,\lambda}^{(r)}(x) \frac{t^n}{n!}.
\end{split}
\end{equation}

From \eqref{eq04} and \eqref{eq05}, we note that
\begin{equation}\label{eq06}
\begin{split}
\phi_{n,\lambda}^{(r)}(x)=\sum_{k=0}^n \begin{Bmatrix}n+r \\ k+r \end{Bmatrix}_{r,\lambda} x^k,\quad (n\geq0), \quad \ \ (\text{see \ \cite{12}}).
\end{split}
\end{equation}

The boson operators $a^\dag$ and $a$ satisfy the following commutation relation:

\begin{equation}\label{eq07}
\begin{split}
[a,a^\dag]=aa^\dag-a^\dag a=1, \quad \ \ (\text{see \ \cite{1, 5, 13}}).
\end{split}
\end{equation}

The number states $|m\rangle, m=0,1,2,\cdots,$ are defined as

\begin{equation}\label{eq08}
\begin{split}
a|m\rangle=\sqrt m|m-1\rangle, \ \ \ a^\dag|m\rangle=\sqrt{m+1}|m+1\rangle, \quad \ \ (\text{see \ \cite{1, 8, 9}}).
\end{split}
\end{equation}

The number operator is defined by

\begin{equation}\label{eq09}
\begin{split}
\hat{n}|k\rangle = k|k\rangle,\quad (k\geq0), \quad \ \ (\text{see \ \cite{5, 8, 9}}).
\end{split}
\end{equation}

By \eqref{eq08} and \eqref{eq09}, we get $\hat{n}=a^\dag a$.

Thus, we note

\begin{equation*}
\begin{split}
[a,\hat{n}]=a\hat{n}-\hat{n}a=a, [\hat{n},a^\dag]=\hat{n}a^\dag-a^\dag \hat{n}=a^\dag.
\end{split}
\end{equation*}

The normal ordering of an integral power of the number operator $\hat{n}=a^\dag a$ in terms of boson operators $a$ and $a^\dag$ can be written in the form:

\begin{equation}\label{eq10}
\begin{split}
(a^\dag a)^k=\sum_{l=0}^k S_2(k,l)(a^\dag)^la^l, \quad \ \ (\text{see \ \cite{1, 5, 8, 13}}),
\end{split}
\end{equation}
where $S_2(k,l)$ are the Stirling numbers of the second kind defined by $x^n=\sum_{k=0}^n S_2(n,k)(x)_k$, \quad (see \ \cite{3, 14}).

The degenerate normal ordering of a degenerate integral power of the number operator $\hat{n}$ in terms of boson operators $a$ and $a^\dag$ is given by

\begin{equation}\label{eq11}
\begin{split}
(a^\dag a)_{k,\lambda}=\sum_{l=0}^k S_{2,\lambda}(k,l)(a^\dag)^la^l, \quad \ \ (\text{see \ \cite{8, 9}}),
\end{split}
\end{equation}
where $S_{2,\lambda}(n,k)$ are the degenerate Stirling numbers of the second kind defined by

\begin{equation}\label{eq12}
\begin{split}
(x)_{n,\lambda}=\sum_{k=0}^nS_{2,\lambda}(n,k)(x)_k, \quad (n\geq0), \quad \ \ (\text{see \ \cite{6}}).
\end{split}
\end{equation}

\section{Degenerate $r$-Bell polynomials arising from degenerate normal ordering}

We recall that the coherent states $|z\rangle, \ z \in \mathbb{C}$, satisfy $a|z\rangle=z|z\rangle, \langle z|z\rangle=1$. Note that $\langle z|a^\dag =\langle z|\bar{z}$.
For the coherent state $|z\rangle,$ we have
\begin{equation}\label{eq13}
\begin{split}
|z\rangle=e^{-\frac{|z|^2}{2}} \sum_{n=0}^\infty \frac{z^n}{\sqrt{n!}}|n\rangle, \quad \ \ (\text{see \ \cite{1, 8, 9, 13}}).
\end{split}
\end{equation}

By \eqref{eq13}, we get

\begin{equation}\label{eq14}
\begin{split}
\langle x|y \rangle=e^{-\frac{1}{2}(|x|^2+|y|^2)+\bar{x}y}, \quad (x,y \in \mathbb{C}).
\end{split}
\end{equation}

It is easy to show that
\begin{equation}\label{eq15}
\begin{split}
\frac{d}{dx}x=1+x\frac{d}{dx}.
\end{split}
\end{equation}

We recall that the standard boson commutation relation $[a,a^\dag]=aa^\dag-a^\dag a=1$ can be considered, in a suitable space of functions $f$, by letting $a=\frac{d}{dx}$ and $a^\dag=x$.

Now we observe that
\begin{equation}\label{eq16}
\begin{split}
\bigg(x\frac{d}{dx}+r\bigg)_{n,\lambda}f(x)=\sum_{k=0}^n \begin{Bmatrix}n+r \\ k+r \end{Bmatrix}_{r,\lambda} x^{k} \bigg(\frac{d}{dx}\bigg)^k f(x),
\end{split}
\end{equation}

and

\begin{equation}\label{eq17}
\begin{split}
\bigg(x\frac{d}{dx}\bigg)_{n,\lambda}x^r f(x)=\sum_{k=0}^n \begin{Bmatrix}n+r \\ k+r \end{Bmatrix}_{r,\lambda} x^{k+r}\bigg(\frac{d}{dx}\bigg)^k f(x),
\end{split}
\end{equation}
where $n,r \in \mathbb{Z}$ with $n,\ r\geq0$.

The equations \eqref{eq16} and \eqref{eq17} can be represented respectively by the degenerate normal orderings of the degenerate $n$-th powers  in \eqref{eq18} and \eqref{eq19} of the number operator $\hat{n}=a^\dag$ in terms of boson operators $a$ and $a^\dag$:
\begin{equation}\label{eq18}
\begin{split}
(\hat{n}+r)_{n,\lambda}=(a^\dag a+r)_{n,\lambda}=\sum_{k=0}^n \begin{Bmatrix}n+r \\ k+r \end{Bmatrix}_{r,\lambda}(a^\dag)^k a^k,
\end{split}
\end{equation}

and

\begin{equation}\label{eq19}
\begin{split}
(\hat{n})_{n,\lambda}(a^\dag)^r=(a^\dag a)_{n,\lambda}(a^\dag)^r=\sum_{k=0}^n \begin{Bmatrix}n+r \\ k+r \end{Bmatrix}_{r,\lambda} (a^\dag)^{k+r} a^k.
\end{split}
\end{equation}

Therefore, from \eqref{eq18} and \eqref{eq19}, we obtain the following theorem.

\begin{theorem}
For $n\geq0$, we have
\begin{equation*}
\begin{split}
&(\hat{n}+r)_{n,\lambda}=(a^\dag a+r)_{n,\lambda}=\sum_{k=0}^n \begin{Bmatrix}n+r \\ k+r \end{Bmatrix}_{r,\lambda}(a^\dag)^k a^k,
\end{split}
\end{equation*}

and, for $n \ge r$, we have

\begin{equation*}
\begin{split}
&(\hat{n})_{n-r,\lambda}(a^\dag)^r a^r=(a^\dag a)_{n-r,\lambda}(a^\dag)^r a^r=\sum_{k=r}^n \begin{Bmatrix}n \\ k \end{Bmatrix}_{r,\lambda} (a^\dag)^{k} a^k.
\end{split}
\end{equation*}
\end{theorem}

Let $m=0,1,2,\cdots.$ Then, by \eqref{eq08}, we get
\begin{equation}\label{eq20}
\begin{split}
(a^\dag a+r)_{n,\lambda}|m\rangle=(\hat{n}+r)_{n,\lambda}|m\rangle = (m+r)_{n,\lambda}|m\rangle,
\end{split}
\end{equation}

and

\begin{equation}\label{eq21}
\begin{split}
(\hat{n}+r)_{n,\lambda}|m\rangle &= (a^\dag a+r)_{n,\lambda}|m\rangle = \sum_{k=0}^n \begin{Bmatrix}n+r \\ k+r \end{Bmatrix}_{r,\lambda} (a^\dag)^ka^k|m\rangle \\
&=\sum_{k=0}^n \begin{Bmatrix}n+r \\ k+r \end{Bmatrix}_{r,\lambda}(m)_k|m\rangle.
\end{split}
\end{equation}

Thus, by \eqref{eq20} and \eqref{eq21}, we get
\begin{equation}\label{eq22}
\begin{split}
(m+r)_{n,\lambda}=\sum_{k=0}^n \begin{Bmatrix}n+r \\ k+r \end{Bmatrix}_{r,\lambda}(m)_k, \quad (n\geq0).
\end{split}
\end{equation}

This is the classical expression for the degenerate $n$-th power of $m+r$ in terms of the falling factorial $(m)_k$. This shows that \eqref{eq03} holds for all nonnegative integers $x=m=0,1,2,\cdots,$ which in turn implies \eqref{eq03} itself holds true.

From \eqref{eq22}, we note that

\begin{equation}\label{eq23}
\begin{split}
\langle z|(\hat{n}+r)_{n,\lambda}|z \rangle&= \langle z|(a^\dag a+r)_{n,\lambda}|z \rangle \\
&=\sum_{k=0}^n \begin{Bmatrix}n+r \\ k+r \end{Bmatrix}_{r,\lambda} \langle z|(a^\dag)^k a^k|z\rangle \\
&=\sum_{k=0}^n \begin{Bmatrix}n+r \\ k+r \end{Bmatrix}_{r,\lambda} (\bar{z})^kz^k \langle z|z\rangle \\
&=\sum_{k=0}^n \begin{Bmatrix}n+r \\ k+r \end{Bmatrix}_{r,\lambda} |z|^{2k}=\phi_{n,\lambda}^{(r)}(|z|^2).
\end{split}
\end{equation}

Therefore, by \eqref{eq23}, we obtain the following theorem.

\begin{theorem}

For $n\geq0$, we have
\begin{equation*}
\begin{split}
\langle z|(\hat{n}+r)_{n,\lambda}|z\rangle= \langle z|(a^\dag a+r)_{n,\lambda}|z \rangle = \phi_{n,\lambda}^{(r)}(|z|^2).
\end{split}
\end{equation*}

\end{theorem}

%8

Let us take $f(t)=\langle \  z \ | \ e_{\lambda}^{r+\hat{n}}(t) \ | \ z \ \rangle$. Then by \eqref{eq23}, we get

\begin{equation}\label{eq24}
\begin{split}
f(t)=\langle \ z \ | \ e_\lambda^{r+\hat{n}}(t) \ | \ z \ \rangle&=\sum_{k=0}^\infty \frac{t^k}{k!}\langle \  z \ | \ (\hat{n}+r)_{k,\lambda} \ | \ z \ \rangle=\sum_{k=0}^\infty \frac{t^k}{k!}\langle \ z \ | \ (a^\dag a+r)_{k,\lambda} \ | \ z \ \rangle\\
&=\sum_{k=0}^\infty \frac{t^k}{k!}\sum_{l=0}^k\begin{Bmatrix} k+r \\ l+r \end{Bmatrix}_{r,\lambda}|z|^{2k}=\sum_{k=0}^\infty \phi_{k,\lambda}^{(r)}(|z|^2)\frac{t^k}{k!}.
\end{split}
\end{equation}

Indeed, the equation \eqref{eq24} says that $f(t)$ is the generating function of the degenerate $r$-Bell polynomials which are considered by Kim-Kim.

Therefore, by \eqref{eq24}, we obtain the following theorem.

\begin{theorem}

The generating function of degenerate $r$-Bell polynomials is given by
\begin{equation*}
\begin{split}
\langle \ z \ | \ e_\lambda^{r+\hat{n}}(t) \ | \ z \ \rangle = \langle \ z \ | \ e_\lambda^{r+a^\dag a}(t) \ | \ z \ \rangle=\sum_{n=0}^\infty \phi_{n,\lambda}^{(r)}(|z|^2)\frac{t^n}{n!}.
\end{split}
\end{equation*}

\end{theorem}

%9
Now, we observe that
\begin{equation}\label{eq25}
\begin{split}
\hat{n}(\hat{n}+r-\lambda)_{k,\lambda}&=a^\dag a(a^\dag a+r-\lambda)_{k,\lambda}\\
&=a^\dag a(a^\dag a+r-\lambda)(a^\dag a+r-2\lambda)\cdots(a^\dag a+r-k\lambda)\\
&=a^\dag(aa^\dag+r-\lambda)a(a^\dag a+r-2\lambda)\cdots (a^\dag a+r-k\lambda)\\
&=a^\dag(a^\dag a+1+r-\lambda)a(a^\dag a+r-2\lambda)\cdots(a^\dag a+r-k\lambda)\\
&=\cdots\\
&=a^\dag(a^\dag a+1+r-\lambda)(a^\dag a+1+r-2\lambda)\cdots(a^\dag a+1+r-k\lambda)a\\
&=a^\dag(a^\dag a+1+r-\lambda)_{k,\lambda}a=a^\dag(\hat{n}+1+r-\lambda)_{k,\lambda}a,\quad (k\geq0).
\end{split}
\end{equation}

From \eqref{eq25}, we note that
\begin{equation}\label{eq26}
\begin{split}
\hat{n}e_\lambda^{\hat{n}+r-\lambda}(t)=e_\lambda^{\hat{n}+r-\lambda}(t)\hat{n}=a^\dag e_\lambda^{aa^\dag+r-\lambda}(t)a=a^\dag e_\lambda^{\hat{n}+1+r-\lambda}(t)a.
\end{split}
\end{equation}

By \eqref{eq24} and \eqref{eq26}, we get
\begin{equation}\label{eq27}
\begin{split}
\frac{\partial f(t)}{\partial t}&=\frac{\partial}{\partial t}\langle \ z \ | \ e_{\lambda}^{a^\dag a+r}(t) \ | \ z \ \rangle=\langle \ z \ | \ (\hat{n}+r)e_{\lambda}^{\hat{n}+r-\lambda}(t) \ | \ z \ \rangle\\
&=\langle \ z \ | \ a^\dag e_{\lambda}^{a^\dag a+r+1-\lambda}(t) a \ | \ z \ \rangle+ r \langle \ z \ | \ e_{\lambda}^{a^\dag a+r-\lambda}(t) \ | \ z \ \rangle\\
&=e_{\lambda}^{1-\lambda}(t)\bar{z}z\langle \ z \ | \ e_{\lambda}^{a^\dag a+r}(t) \ | \ z \ \rangle+ r\langle \ z \ | \ e_{\lambda}^{a^\dag a+r}(t) \ | \ z \ \rangle e_\lambda^{-\lambda}(t)\\
&=e_\lambda^{1-\lambda}(t)|z|^2f(t)+re_\lambda^{-\lambda}(t)f(t).
\end{split}
\end{equation}

Therefore, by \eqref{eq27}, we obtain the following theorem.

%10
\begin{theorem}

Let $f(t)=\langle \ z \ | \ e_{\lambda}^{r+\hat{n}}(t) \ | \ z  \ \rangle$. Then we have
\begin{equation*}
\begin{split}
\frac{f'(t)}{f(t)}=|z|^2e_\lambda^{1-\lambda}(t)+re_\lambda^{-\lambda}(t).
\end{split}
\end{equation*}

\end{theorem}

Note that $f(0)=1$. From Theorem 4, we have
\begin{equation}\label{eq28}
\begin{split}
\log f(t)&=\int_0^t\frac{f'(t)}{f(t)}dt=|z|^2\int_0^t e_\lambda^{1-\lambda}(t)dt+r\int_0^te_\lambda^{-\lambda}(t)dt\\
&=|z|^2(e_\lambda(t)-1)+r\frac{1}{\lambda}\log(1+\lambda t).
\end{split}
\end{equation}

Thus, by \eqref{eq28}, we get
\begin{equation}\label{eq29}
\begin{split}
f(t)=e^{|z|^2(e_\lambda(t)-1)+r\frac{1}{\lambda}\log(1+\lambda t)}=e_\lambda^{r}(t)e^{|z|^2(e_\lambda(t)-1)}.
\end{split}
\end{equation}

Therefore, by \eqref{eq29}, we obtain the following theorem.

\begin{theorem}

Let $f(t)=\langle \ z \ | \ e_{\lambda}^{r+\hat{n}}(t) \ | \ z \ \rangle=\langle \ z \ | \ e_\lambda^{r+a^\dag a} \ | \ z \ \rangle$. Then we have
\begin{equation*}
\begin{split}
f(t)=e_\lambda^r(t)e^{|z|^2(e_\lambda(t)-1)}.
\end{split}
\end{equation*}

\end{theorem}

From Theorems 3 and 5, we have
%11
\begin{equation}\label{eq30}
\begin{split}
e_\lambda^r(t)e^{|z|^2(e_\lambda(t)-1)}=\sum_{n=0}^\infty\phi_{n,\lambda}^{(r)}(|z|^2)\frac{t^n}{n!}.
\end{split}
\end{equation}

That is, by \eqref{eq29}, we get
\begin{equation}\label{eq31}
\begin{split}
\sum_{n=0}^\infty \langle \ z \ | \ (r+\hat{n})_{n,\lambda} \ | \ z \ \rangle\frac{t^n}{n!}&=\langle \ z \ | \ e_\lambda^{r+\hat{n}}(t) \ | \ z \ \rangle =f(t)\\
&=\sum_{k=0}^\infty |z|^{2k}\frac{1}{k!}(e_\lambda(t)-1)^ke_\lambda^r(t)\\
&=\sum_{k=0}^\infty |z|^{2k}\sum_{n=k}^\infty \begin{Bmatrix} n+r \\ k+r \end{Bmatrix}_{r,\lambda}\frac{t^n}{n!}\\
&=\sum_{n=0}^\infty \bigg(\sum_{k=0}^n |z|^{2k}\begin{Bmatrix} n+r \\ k+r \end{Bmatrix}_{r,\lambda}\bigg)\frac{t^n}{n!}.
\end{split}
\end{equation}

Comparing the coefficients on both sides of \eqref{eq31}, we have
\begin{equation}\label{eq32}
\begin{split}
\langle \ z \ | \ (r+\hat{n})_{n,\lambda} \ | \ z \ \rangle=\langle \ z \ | \ (r+a^\dag a)_{n,\lambda} \ | \ z \ \rangle=\sum_{k=0}^n\begin{Bmatrix} n+r \\ k+r \end{Bmatrix}_{r,\lambda}|z|^{2k}.
\end{split}
\end{equation}

By Theorems 3 and 5, we get
\begin{equation}\label{eq33}
\begin{split}
f(t)=\langle \ z \ | \ e_\lambda^{r+a^\dag a}(t) \ | \ z \ \rangle=\sum_{n=0}^\infty \phi_{n,\lambda}^{(r)}(|z|^2)\frac{t^n}{n!}.
\end{split}
\end{equation}

Differentiating \eqref{eq33} with respect to $t$, we have
\begin{equation}\label{eq34}
\begin{split}
\frac{\partial f(t)}{\partial t}=\sum_{k=1}^\infty \frac{t^{k-1}}{(k-1)!}\phi_{k,\lambda}^{(r)}(|z|^2)=\sum_{k=0}^\infty \phi_{k+1,\lambda}^{(r)}(|z|^2)\frac{t^k}{k!}.
\end{split}
\end{equation}

On the other hand, by \eqref{eq27}, we get
%12
\begin{equation}\label{eq35}
\begin{split}
\frac{\partial f(t)}{\partial t}&=e_\lambda^{1-\lambda}(t)|z|^2f(t)+re_\lambda^{-\lambda}(t)f(t)\\
&=e_\lambda^{1-\lambda}(t)|z|^2\sum_{k=0}^\infty \frac{t^k}{k!}\phi_{k,\lambda}^{(r)}(|z|^2)+re_\lambda^{-\lambda}(t)\sum_{k=0}^\infty\phi_{k,\lambda}^{(r)}(|z|^2)\frac{t^k}{k!}\\
&=|z|^2\sum_{l=0}^\infty (1-\lambda)_{l,\lambda}\frac{t^l}{l!}\sum_{k=0}^\infty \frac{t^k}{k!}\phi_{k,\lambda}^{(r)}(|z|^2)+r\sum_{l=0}^\infty(-\lambda)_{l,\lambda}\frac{t^l}{l!}\sum_{k=0}^\infty\phi_{k,\lambda}^{(r)}(|z|^2)\frac{t^k}{k!}\\
&=\sum_{n=0}^\infty\bigg(|z|^2\sum_{k=0}^n\binom{n}{k}(1-\lambda)_{n-k,\lambda}\phi_{k,\lambda}^{(r)}(|z|^2)+r\sum_{k=0}^n\binom{n}{k}(-\lambda)_{n-k,\lambda}\phi_{k,\lambda}^{(r)}(|z|^2)\bigg)\frac{t^n}{n!}.
\end{split}
\end{equation}

Therefore, by \eqref{eq34} and \eqref{eq35}, we obtain the following theorem.
\begin{theorem}

For $n\geq0$, we have
\begin{equation*}
\begin{split}
\phi_{n+1,\lambda}^{(r)}(|z|^2)=\sum_{k=0}^n\binom{n}{k}\phi_{k,\lambda}^{(r)}(|z|^2)(|z|^2(1-\lambda)_{n-k,\lambda}+r(-\lambda)_{n-k,\lambda}).
\end{split}
\end{equation*}

\end{theorem}

When $|z|=1$ in Theorem 6, we have
\begin{equation*}
\begin{split}
\phi_{n+1,\lambda}^{(r)}=\sum_{k=0}^n\binom{n}{k}\phi_{k,\lambda}^{(r)}((1-\lambda)_{n-k,\lambda}+r(-\lambda)_{n-k,\lambda}),
\end{split}
\end{equation*}

where $\phi_{n,\lambda}^{(r)}(1)=\phi_{n,\lambda}^{(r)}$ are called the degenerate $r$-Bell numbers.

%13

Evaluating the left hand side of Theorem 2 by using the representation of the coherent state in terms of the number operator in \eqref{eq13}, we have
\begin{equation}\label{eq36}
\begin{split}
\langle \ z \ | \ (a^\dag a+r)_{k,\lambda} \ | \ z \ \rangle&=e^{-\frac{|z|^2}{2}}e^{-\frac{|z|^2}{2}}\sum_{m,n=0}^\infty \frac{|\bar{z}|^mz^n}{\sqrt{m!}\sqrt{n!}}(n+r)_{k,\lambda}\langle \ m \ | \ n \ \rangle\\
& = e^{-|z|^2}\sum_{n=0}^\infty \frac{|z|^{2n}}{n!}(n+r)_{k,\lambda}.
\end{split}
\end{equation}

By Theorem 2 and \eqref{eq36}, we get

\begin{equation}\label{eq37}
\begin{split}
\phi_{k,\lambda}^{(r)}(|z|^2)&=\langle \ z \ | \ (a^\dag a+r)_{k,\lambda}\ | \ z \ \rangle =\sum_{l=0}^k |z|^{2l}\begin{Bmatrix} k+r \\ l+r \end{Bmatrix}_{r,\lambda}\\
&=e^{-|z|^2}\sum_{n=0}^\infty \frac{|z|^{2n}}{n!}(n+r)_{k,\lambda},\quad \ \ (k\geq0).
\end{split}
\end{equation}
In particular, by  letting $|z|=1$, we get the Dobinski-like formula.

\begin{theorem}
The following Dobinski-like formual holds true.
\begin{equation*}
\begin{split}
\phi_{k,\lambda}^{(r)}=\frac{1}{e}\sum_{n=0}^\infty\frac{1}{n!}(n+r)_{k,\lambda},\quad \ \ (k, r \in\mathbb{N}\cup\{0\}).
\end{split}
\end{equation*}
\end{theorem}

\medskip
\section{Conclusion}

In recent years, various degenerate versions of many special numbers and polynomials have been explored  by employing diverse tools such as generating functions, combinatorial methods, umbral calculus, $p$-adic analysis, differential equations, probability theory, operator theory and analytic number theory. \par
In this paper, we studied the degenerate $r$-Bell polynomials and numbers via boson operators in quantum physics. To do that, we first recalled the degenerate normal ordering of a degenerate integral power of the number operator in terms of boson operators. Here the degenerate $r$-Stirling numbers of the second kind appear as the coefficients. From the degenerate normal ordering, we derived two expressions for the generating function of the degenerate $r$-Bell polynomials in $|z|^{2}$, and a recurrence relation and Dobinski-like formula for the degenerate $r$-Bell numbers. \par
It is one of our future research projects to continue to study many different versions of some special numbers and polynomials and to find their applications in physics, science and engineering as well as in mathematics.

\bigskip

\noindent{\bf{Acknowledgments}} \\
%The author would like to thank the referees for the detailed and valuable comments that helped
%improve the original manuscript in its present form.% Also,
%%
The authors thank Jangjeon Institute for Mathematical Science for the support of this research.

\vspace{0.1in}

\noindent{\bf {Availability of data and material}} \\
Not applicable.

\vspace{0.1in}

\noindent{\bf{Funding}} \\
This work  supported by the Basic Science Research Program, the National
Research Foundation of Korea, (NRF-2021R1F1A1050151).
\vspace{0.1in}

\noindent{\bf{Ethics approval and consent to participate}} \\
All authors declare that there is no ethical problem in the production of this paper.

\vspace{0.1in}

\noindent{\bf {Competing interests}} \\
All authors declare no conflict of interest.

\vspace{0.1in}

%\noindent{\bf{Consent for publication}} \\
%All authors want to publish this paper in this journal.
%
%\vspace{0.1in}
%
\noindent{\bf{Author' Contributions}}\\
TK and DSK conceived of the framework and structured the whole
paper; DSK and TK wrote the whole paper. HKK checked the results
of the paper and typed the paper: HKK paid for the charge of
publication fee of this paper. All authors read and agreed with
the published version of the manuscript.
%\vspace{0.1in}
%
%\noindent{\bf{Author details}}\\
%Department of Mathematics Education, Daegu Catholic University, Gyeongsan 38430, Republic of Korea.

%\vspace{0.1in}
\


\bigskip
\begin{thebibliography}{9}


\bibitem{1}
Blasiak, P. Combinatorics of Boson normal ordering and some applications. Ph D Thesis, University of Paris, Paris, (2005). http://arxiv.org/pdf/quant-ph/0507206.pdf

\bibitem{2}
Carlitz, L. Degenerate Stirling, Bernoulli and Eulerian numbers. Utilitas Math. 15 (1979), 51-88.


\bibitem{3}
Comtet, L. Advanced combinatorics. The art of finite and infinite expansions. Revised and enlarged edition. D. Reidel Publishing Co., Dordrecht, 1974. xi+343 pp. ISBN: 90-277-0441-4 05-02


\bibitem{4}
El-Desouky, B. S. Multiparameter non-central Stirling numbers. Fibonacci Quart. 32 (1994), no. 3, 218-225.


\bibitem{5}
Katriel, J. Sirling number identities: inter consistency of $q$-analogues. J. Phys. A 31 (1998), no. 15, 3559-3572.

\bibitem{6}
Kim, D. S.;  Kim, T.  A note on a new  type of degenerate Bernoulli numbers, Russ. J. Math. Phys. 27 (2020), no. 2, 227-235.

\bibitem{7}
Kim, T.; Kim, D. S. On some degenerate differential and degenerate difference operators. Russ. J. Math. Phys. 29 (2022), no. 1, 37-46.


\bibitem{8}
Kim, T.; Kim, D. S. Degenerate $r$-Whitney numbers and degenerate $r$-Dowling polynomials via boson operators. Adv. in Appl. Math. 140 (2022), Paper No. 102394.


\bibitem{9}
Kim, T.; Kim, D. S.; Kim, H. K. Normal ordering of degenerate integral powers of number operator and its applications. Appl. Math. Sci. Eng. 30 (2022), no. 1, 440-447.


\bibitem{10}
Kim, T.; Kim, D. S.; Kwon, J.; Lee, H. Some identities involving degenerate $r$-Stirling numbers. Proc. Jangjeon Math. Soc. 25 (2022), no. 2, 245-252.


\bibitem{11}
Kim, T.; Kim, D. S.; Lee, H.; Park, J.-W. A note on degenerate $r$-Stirling numbers. J. Inequal. Appl. 2020, Paper No. 225, 12 pp.


\bibitem{12}
Kim, T.; Yao, Y.; Kim, D. S.; Jang. G.-W. Degenerate $r$-Stirling numbers and $r$-Bell polynomials. Russ. J. Math. Phys. 25 (2018), no. 1, 44-58.



\bibitem{13}
Perelomov, A. Generalized coherent states and their applications. Texts and Monographs in Physics. Springer Verlag, Berlin, 1986. xii+320 pp.


\bibitem{14}
Roman, S. The umbral Calculus. Pure and Applied Mathematics, 111. Academic Press, Inc. [Harcourt Brace Jovanovich, Publishers], New York, 1984. x+193 pp.


\end{thebibliography}
\end{document}